\documentclass[12pt,amsmath]{article}
\usepackage{amssymb}
\usepackage{amscd}
\usepackage{times}
\begin{document}
\newtheorem{theorem}{\Large\bf\sf Theorem}
\newtheorem{thm}[theorem]{\Large\bf\sf Theorem}
\newtheorem{proposition}[theorem]{\Large\bf\sf Proposition}
\newtheorem{corollary}[theorem]{\Large\bf\sf Corollary}
\newtheorem{lemma}[theorem]{\Large\bf\sf Lemma}
\newtheorem{lem}[theorem]{\Large\bf\sf Lemma}
\newtheorem{remark}[theorem]{\Large\bf\sf Remark}

\def\K{\mathop{\mathsf K}\nolimits}
\def\Id{\mathop{\operatorname{\rm Id}\nolimits}}
\def\ch{\mathop{\operatorname{\rm ch}\nolimits}}
\def\Ind{\mathop{\operatorname{\rm Ind}\nolimits}}
\def\Mat{\mathop{\mathsf {Mat}}\nolimits}
\def\GL{\mathop{\mathsf {GL}}\nolimits}
\def\SL{\mathop{\operatorname{\rm SL}\nolimits}}
\def\HE{\mathop{\operatorname{\rm HE}\nolimits}}
\def\HP{\mathop{\mathsf {HP}}\nolimits}
\def\HC{\mathop{\mathsf {HC}}\nolimits}
\def\HH{\mathop{\mathsf {HH}}\nolimits}
\def\Hom{\mathop{\operatorname{\rm Hom}\nolimits}}
\def\Tot{\mathop{\operatorname{\rm Tot}\nolimits}}
\def\ad{\mathop{\operatorname{\rm ad}\nolimits}}
\def\Ad{\mathop{\operatorname{\rm Ad}\nolimits}}
\def\Aff{\mathop{\operatorname{\rm Aff}\nolimits}}
\def\aff{\mathop{\operatorname{\rm aff}\nolimits}}
\def\Aut{\mathop{\operatorname{\rm Aut}\nolimits}}
\def\Lie{\mathop{\operatorname{\rm Lie}\nolimits}}

\def\Ad{\mathop{\mathrm {Ad}}\nolimits}
\def\ad{\mathop{\mathrm {ad}}\nolimits}
\def\Lie{\mathop{\operatorname {Lie}}\nolimits}
\def\Aff{\mathop{\mathsf {Aff}}\nolimits}
\def\aff{\mathop{\mathsf {aff}}\nolimits}
\def\SL{\mathop{\mathsf {SL}}\nolimits}
\def\Id{\mathop{\mathsf {Id}}\nolimits}
\def\Ad{\mathop{\mathsf {Ad}}\nolimits}
\def\ad{\mathop{\mathsf {ad}}\nolimits}
\def\Ind{\mathop{\mathsf {Ind}}\nolimits}
\def\Id{\mathop{\mathsf {Id}}\nolimits}
\def\Aut{\mathop{\mathsf {Aut}}\nolimits}
\def\End{\mathop{\mathsf {End}}\nolimits}
\def\sgn{\mathop{\mathsf {sgn}}\nolimits}
\def\Lie{\mathop{\mathsf {Lie}}\nolimits}
\def\Ln{\mathop{\mathsf {Ln}}\nolimits}
\def\tr{\mathop{\mathsf {tr}}\nolimits}
\def\Aff{\mathop{\mathsf {Aff}}\nolimits}
\def\aff{\mathop{\mathsf {aff}}\nolimits}
\def\SL{\mathop{\mathsf {SL}}\nolimits}
\def\Id{\mathop{\mathsf {Id}}\nolimits}
\def\Ad{\mathop{\mathsf {Ad}}\nolimits}
\def\ad{\mathop{\mathsf {ad}}\nolimits}
\def\Ind{\mathop{\mathsf {Ind}}\nolimits}
\def\Id{\mathop{\mathsf {Id}}\nolimits}
\def\Aut{\mathop{\mathsf {Aut}}\nolimits}
\def\End{\mathop{\mathsf {End}}\nolimits}
\def\sgn{\mathop{\mathsf {sgn}}\nolimits}

\title{\bf\sf Riemann-Roch Theorem and Index Theorem in Non-commutative Geometry
\footnote{\large The text was presented at the ICM2002-Satellite Conference ``Abstract and Applied Analysis 2002", August 13-17, 2002, Hanoi, the author expresses his sincere and deep thanks to the organizers, and especially Professor Dr. DSc. Nguyen Minh Chuong for invitation}}
\author{\bf\sf Do Ngoc Diep}
\date{\bf Institute of Mathematics\\ National Centre for Science and Technology of Vietnam}
\maketitle
\section{\bf The classical Riemann-Roch Theorem}\Large\bf\sf  The classical Riemann-Roch Theorem is well-known in the function theory as [M.F. Atiyah and Hirzebruch, {\it Riemann-Roch Theorem for differential manifolds}, Bull. Amer. Math. Soc. 1959, Vol. 65, 276-281]
\begin{theorem}[\Large\bf\sf The Riemann-Roch Theorem]
$$r(-D) -i(D) = d(D) - g + 1,$$
where $D$ is a fixed divisor of degree $d(D)$ on a Riemann surface $X$ of genus $g$, $r(-D)$ is the dimension of the space of meromorphic functions of divisor $\geq -D$ on $X$, $i(D)$ the dimension of the space of meromorphic 1-forms of divisor $\geq D$ on $X$.
\end{theorem}

This theorem can be considered as computing the Euler characteristics of the sheaf of germs of holomorphic sections of the holomorphic bundle, defined by the divisor $D$, over $X$.
It plays also an important role in classical algebraic geometry.

\begin{theorem} \Large\bf\sf
Let $X$ be a nonsingular complex projective algebraic variety, $c$ its first
Chern class, $\xi$ a holomorphic bundle over $X$. Then the value on $[X]$ of the
cohomological class $$e^{c/2}. ch\xi . A^{-1}(p_1(X),p_2(X),\dots)$$ equal to
the Euler characteristic of the sheaf of holomorphic sections of the bundle
$\xi$. \end{theorem}

\section{\bf The Riemann-Roch Theorem in Algebraic Topology} 
In algebraic topology the Riemann-Roch theorem appeared as some measure of noncommutativity of some diagrams relating two generalized homology theories.
Let us remind the most general setting of the Riemann-Roch Theorem.
\par
Consider two generalized (co)homology theories $k(X)$ and $h(X)$ and $\tau: h \to k$ be a multiplicative map, sending $1\in h^0(pt)$ to $1\in k^0(pt)$. Consider a vector bundle $\xi$, oriented with respect to the both theories, over the base $X$. Denote $T(\xi)$ the {\it Thom space} of $\xi$, i.e. the quotient of the corresponding disk bundle $D(\xi)$ modulo its boundary $Sph(\xi)$.
Let us consider an $h$-oriented vector bundle $\xi = (E,X,V,p)$. Denote $E'$ the complement to the zero section of $E$. It is easy to see that $\tilde{h}^r(T(\xi)) = h^r(E,E')$.
\begin{theorem}[\Large\bf\sf Thom Isomorphism]    \Large\bf\sf
The Thom homomorphism $$ \CD t: h^q(X) =h^q(E) @>.u\in h^n(E,E') >> h^{q+n}(E,E') = \tilde{h}^{q+n}(T(\xi)) \endCD$$
is an isomorphism.
\end{theorem}
Following the Thom isomorphism theorem, there are Thom isomorphisms $t_h^\xi: \tilde{h}^*(T(\xi)) \to h^*(X)$ and $t_k^\xi: \tilde{k}^*(T(\xi)) \to \tilde{k}^*(X)$. The Todd class is defined as $\mathcal T_\tau(\xi) := (t_k^\xi)^{-1} \circ \tau \circ t_h^\xi (1)$. The most general Riemann-Roch Theorem states: 

\begin{theorem}  \Large\bf\sf
For every $\alpha\in h^*(X)$, one has $$(t^\xi_k)^{-1} \circ \tau \circ t_h^\xi(\alpha) = \tau(\alpha) \mathcal T_\tau(\xi)$$
\end{theorem} 
The Todd class is therefore some noncommutativity measure of the diagram
$$\CD
\tilde{h}^*(T(\xi)) @>\tau>> \tilde{k}^*(T(\xi))\\
@AA t_h^\xi A   @AA t_k^\xi A\\
h^*(X) @>\tau >> k^*(X) \endCD$$
Example: $h=k=H^*(.,\mathbb Z_2)$, $\tau = Sq = 1+Sq^1+Sq^2+\dots$, then $\mathcal T_\tau^\xi = w(\xi) = 1+w_1(\xi) + w_2(\xi) + \dots$ - the full Stiefel-Whitney class.

\section{\bf The Riemann-Roch Theorem and the Index Theorem of pseudo-differential operators}
One of the consequences of the Riemann-Roch Theorem is the fact that the index of the Dirac operator on $X$ is exactly the Euler characteristics of $X$.
Let us review the classical results from algebraic topology and topology of pseudo-differential operators.
With the Riemann-Roch Theorem it is convenient to define the {\it direct image map} $f_!$ as the special case of the composition map
$$\CD 
h^*(X) @>t_h >> \tilde{h}^*(T(\nu(X))) @>g^* >> \tilde{h}^*(T(\nu(Y))) @>t^{-1}_h >> h^*(Y),
\endCD$$ where $h$ is the ordinary co-homology, $X$ and $Y$ are the oriented manifolds and $g=f: X \to Y$.

\begin{theorem}\Large\bf\sf  Let $X$ and $Y$ be two closed manifold, oriented with respect to the both generalized (co)homology theories $h$ and $k$ and $f: X \to Y$ a continuous map. Then the Todd classes $\mathcal T(X)$ and $\mathcal T(Y)$ measure noncommutativity of the diagram
$$\CD
h^*(X) @>\tau >> k^*(X)\\
@VV f_!V    @VV f_! V\\
h^*(Y) @>\tau >> h^*(Y)\endCD$$ and more precisely, $f_!(\tau(x)\mathcal T(X)) =\tau(f_!(x)).\mathcal T(Y)$.
\end{theorem}
The index of an arbitrary pseudo-differential operator $D$ is reduced to the index of Dirac operator $d+\delta$.
\begin{theorem}[\Large\bf\sf Atiyah-Singer-Hirzebruch Index Theorem]   \Large\bf\sf
$$index D = \langle (ch\;D\circ \mathcal T^{-1}((C_\tau(X))),[X]\rangle,$$ where $\mathcal T(\mathbb C_\tau(X))$ is the Todd class of complexified tangent bundle, $\mathcal T^{-1}(\mathbb C_\tau(X)) = U(p_1(x),p_2(x),\dots)$, and $U(e_1(x^2),e_2(x^2),\dots) = - \prod_i \frac{x^2}{(1-e^{-x_i})(1-e^{-x_i}}$.
\end{theorem}
\par
As a consequence of the previous theorem we have the following result.
Let $X$ denote a $2n$-dimensional oriented closed smooth manifold with spin structure, i.e. a fixed Hermitian structure on fibers, smoothly depending of points on the base $X$. Denote $\Omega^k(X)$ the space of alternating differential $k$-forms on $X$, and $d: \Omega^k(X) \to \Omega^{k+1}(X)$ the exterior differential,
${}^*$ the Hodge star operator and $\delta ={}^*d{}^*: \Omega^{k-1}(X) \to \Omega^k(X)$ the adjoint to $d$ operator.
The Dirac operator is $d+\delta$ is a first order elliptic operator and its index is just equal to the Euler characteristic $\chi(X)$ of the manifold.
$$\chi(X) := \sum_r (-1)^r \dim H^r(X,\eta).$$
\begin{theorem} \Large\bf\sf
$$ind (d+\delta) = \chi(X) $$
\end{theorem}

\section{\bf Riemann-Roch Theorem in non-commutative geometry}
Let us consider an arbitrary algebra $A$ over the ground field of complex numbers $\mathbb C$ and $G$ a locally compact group and denote $dg$ the left-invariant Haar measure on $G$. The space $C^\infty_c(G)$ of all continuous functions with compact support on $G$ with values in $A$ under ordinary convolution
$$(f*g) (x) := \int_G f(y)g(y^{-1}x)dy, $$
involution
$$f^*(x) = \overline{f(x^{-1})}$$
and norm
$$\Vert f \Vert = \sup_{x\in G}\vert f(x) \vert$$
form a the so called cross-product $A \rtimes G$ of $A$ and $G$.
\begin{theorem}[\Large\bf\sf Connes-Thom Isomorphism]    \Large\bf\sf
$$t_{\K}: \K_{*+1}(A \rtimes \mathbb R) \cong \K_{*}(A),$$ $$ \quad t_{\HC}: \HC_{*+1}(A \rtimes \mathbb R) \cong \HC_{*}(A),$$ $$ t_{\HP}: \HP_{*+1}(A \rtimes \mathbb R) \cong \HP_*(A).$$
\end{theorem}
A consequence of this theorem is the existence of some noncommutative Todd class
\begin{theorem}\Large\bf\sf  There exists some Todd class
$$(\mathcal T_\tau =(t_K)^{-1} \circ \tau \circ t_{\HP} (1)$$ which measures the noncommutativity of the diagram
$$\CD
\K_{*+1}(A \rtimes \mathbb R) @>\tau >> \HP_{*+1}(A\rtimes \mathbb R)\\
@VV t_K^{\mathbb R} V      @VV t_{\HP} V\\
\K_*(A) @>\tau >> \HP_*(A) \endCD$$
\end{theorem}

\section{Deformation quantization and periodic cyclic homology}
Deformation quantization gives us some noncommutative algebras which are deformation of the classical algebras of holomorphic functions on $X$. For an arbitrary noncommutative algebra $A$ there are at least two generalized homology
theories: the K-theory and periodic cyclic homology. The Connes-Thom isomorphism gives us a possibility to compare the two theories. There appeared some Todd class as the measure of noncommutativity.
Let us review some results of P.Bressler, R. Nest and B. Tsygan [alg-geom/9705014v2 3 Jun 1997].
\par
Deformation quantization of a manifold $M$ is a formal one parameter deformation of the structure sheaf $\mathcal O_M$, i.e. a sheaf of algebras $\mathbb A_M^h$ flat over $\mathbb C[[h]]$ together with an isomorphism of algebras $\mathbb A^h_M \otimes_{\mathbb C[[h]]} \mathbb C \to \mathcal O_M$.
The formula $$\{f,g\} = \frac{1}{h}[\hat{f},\hat{g}] + h.\mathbb A^h_M,$$ where $f$ and $g$ are two local sections of $\mathcal O_M$ and $\hat{f}$, $\hat{g}$ are their respective lifts in $\mathbb A^h_M$, defines a Poisson structure associated to the deformation quantization $\mathbb A^h_M$.
\par
It is well-known that all symplectic deformation quantization of $M$ of dimension $\dim_{\mathbb C} M = 2d$ are locally isomorphic to the standard deformation quantization of $\mathbb C^{2d}$, i.e. in a neighborhood $U'$ of the origin in $\mathbb C^{2d}$ there is an isomorphism 
$$\mathbb A^h_{\mathbb C^{2d}}(U') = \mathcal O_{\mathbb C^{2d}}(U')[[h]] \cong \mathbb A^h_M(U)$$ of algebras over $\mathbb C[[h]]$, continuous in the $h$-adic topology, where the product on $\mathbb A^h_{\mathbb C^{2d}}(U')$ is given in coordinates $x_1,\dots,x_d,\xi_1,\dots,\xi_d$ on $\mathbb C^{2d}$ by the standard Weyl product
$(f*g)(x,\xi) =$ $$ \exp{\left( \frac{\sqrt{-1}h}{2}\sum_{i=1}^d \left(\frac{\partial}{\partial \xi_i}\frac{\partial}{\partial y_i} - \frac{\partial}{\partial\eta_i}\frac{\partial}{\partial x_i}\right)\right)}f(x,\xi)g(y,\eta)\vert_{x=y\atop\xi=\eta} $$
\begin{theorem}[\Large\bf\sf RRT for periodic cyclic cycles]\Large\bf\sf  The diagram
$$\CD
CC^{per}_*(\mathbb A^h_M) @>\sigma>> CC^{per}_*(\mathcal O_M)\\
@Vi VV @VV{\mu_{\mathcal O} \cup \hat{A}(TM) \cup e^\theta} V\\
CC^{per}(\mathbb A^h_M)[h^{-1}] @>\mu^h_{\mathbb A} >> \prod_{p=-\infty}^\infty \mathbb C_M[h^{-1},h][-2p] 
\endCD$$ is commutative.
\end{theorem} 
For $\mathcal D$- and $\mathcal E$-modules of the ring of pseudo-differential operators, take $M=T^*X$ for a complex manifold $X$, and $\mathbb A^h_{T^*X}$ is the deformation quantization with the characteristic class $\theta= \frac{1}{2}\pi^*c_1(X)$, then $\hat{A}(TM)\cup e^\theta = \pi^*\mathcal T(TM)$.
After the use of Gelfand-Fuch cohomology the computation become available.

\section{Noncommutative Chern characters of some quantum algebras}
Let us demonstrate the indicated scheme to some concrete cases of quantum half-planes and quantum punctured complex planes.
We indicate some results, computed by Do Ngoc Diep and Nguyen Viet Hai [Contributions to Algebra and Geometry, Vol. 42, No 2] and Do Ngoc Diep and Aderemi O. Kuku [arXiv.org/math.QA/0109042].

{\it Canonical coordinates on the upper half-planes.}
Recall that the Lie algebra $\mathfrak g = \aff(\mathbf  R)$ of affine transformations of the real straight line is described as follows, see for example \cite{diep2}: The Lie group $\Aff(\mathbf R)$ of affine transformations of type $$x \in \mathbf R \mapsto ax + b, \mbox{ for some parameters }a, b \in \mathbf R, a \ne 0.$$ It is well-known that this group $\Aff(\mathbf R)$ is a two dimensional Lie group which is isomorphic to the group of matrices
$$\Aff(\mathbf R) \cong \{\left (\begin{array}{cc} a & b \\ 0 & 1 \end{array} \right) \vert a,b \in \mathbf R , a \ne 0 \}.$$ We consider its connected component $$G= \Aff_0(\mathbf R)= \{\left (\begin{array}{cc} a & b \\ 0 & 1 \end{array} \right) \vert a,b \in \mathbf R, a > 0 \}$$ of identity element. Its Lie algebra is
$$\mathfrak g = \aff(\mathbf R) \cong  \{\left (\begin{array}{cc} \alpha & \beta \\ 0 & 0 \end{array} \right) \vert \alpha, \beta  \in \mathbf R \}$$  admits a basis of two generators $X, Y$ with the only nonzero Lie bracket $[X,Y] = Y$, i.e.
$$\mathfrak g = \aff(\mathbf R) \cong \{ \alpha X + \beta Y \vert [X,Y] = Y, \alpha, \beta \in \mathbf R \}.$$
The co-adjoint action of $G$ on $\mathfrak g^*$ is given (see e.g. \cite{arnalcortet2}, \cite{kirillov1}) by $$\langle K(g)F, Z \rangle = \langle F, \Ad(g^{-1})Z \rangle, \forall F \in \mathfrak g^*, g \in G \mbox{ and } Z \in \mathfrak g.$$ Denote the co-adjoint orbit of $G$ in $\mathfrak g$, passing through $F$ by
$$\Omega_F = K(G)F :=  \{K(g)F \vert F \in G \}.$$ Because the group $G = \Aff_0(\mathbf R)$ is exponential (see \cite{diep2}), for $F \in \mathfrak g^* = \aff(\mathbf R)^*$, we have
$$\Omega_F = \{ K(\exp(U)F | U \in \aff(\mathbf R) \}.$$
It is easy to see that
$$\langle K(\exp U)F, Z \rangle = \langle F, \exp(-\ad_U)Z \rangle.$$ It is easy therefore to see that
$$K(\exp U)F = \langle F, \exp(-\ad_U)X\rangle X^*+\langle F, \exp(-\ad_U)Y\rangle Y^*.$$
For a general element $U = \alpha X + \beta Y \in \mathfrak g$, we have
$$\exp(-\ad_U) = \sum_{n=0}^\infty \frac{1}{n!}\left(\begin{array}{cc}0 & 0 \\ \beta & -\alpha \end{array}\right)^n = \left( \begin{array}{cc} 1 & 0 \\ L & e^{-\alpha} \end{array} \right),$$ where $L = \alpha + \beta + \frac{\alpha}{\beta}(1-e^\beta)$. This means that
$$K(\exp U)F = (\lambda + \mu L) X^* + (\mu e^{\-\alpha})Y^*. $$ From this formula one deduces  \cite{diep2} the following description of all co-adjoint orbits of $G$ in $\mathfrak g^*$:
\begin{itemize}
\item If $\mu = 0$, each point $(x=\lambda , y =0)$ on the abscissa ordinate corresponds to a 0-dimensional co-adjoint orbit $$\Omega_\lambda = \{\lambda X^* \}, \quad \lambda \in \mathbf R .$$
\item For $\mu \ne 0$, there are two 2-dimensional co-adjoint orbits: the upper half-plane $\{(\lambda , \mu) \quad\vert\quad \lambda ,\mu\in \mathbf R , \mu > 0 \}$ corresponds to the co-adjoint orbit
\begin{equation} \Omega_{+} := \{ F = (\lambda + \mu L)X^* + (\mu e^{-\alpha})Y^* \quad \vert \quad \mu > 0 \}, \end{equation}
and the lower half-plane $\{(\lambda , \mu) \quad\vert\quad \lambda ,\mu\in \mathbf R , \mu < 0\}$ corresponds to the co-adjoint orbit
\begin{equation} \Omega_{-} := \{ F = (\lambda + \mu L)X^* + (\mu e^{-\alpha})Y^* \quad \vert \quad \mu < 0 \}. \end{equation}
\end{itemize}

Denote by $\psi$ the indicated symplectomorphism from $\mathbf R^2$ onto $\Omega_+$
$$(p,q) \in \mathbf R^2 \mapsto \psi(p,q):= (p,e^q) \in \Omega_+$$
\begin{proposition}   \Large\bf\sf
1. Hamiltonian function $f_Z = \tilde{Z}$ in canonical coordinates $(p,q)$ of the orbit $\Omega_+$ is of the form $$\tilde{Z}\circ\psi(p,q) = \alpha p + \beta e^q, \mbox{ if  } Z = \left(\begin{array}{cc} \alpha & \beta \\ 0 & 0 \end{array} \right).$$

2. In the canonical coordinates $(p,q)$ of the orbit $\Omega_+$, the Kirillov form $\omega_{Y^*}$ is just the standard form $\omega = dp \wedge dq$.
\end{proposition}

{\it Computation of generators $\hat{\ell}_Z$}
Let us denote by $\Lambda$ the 2-tensor associated with the Kirillov standard form $\omega = dp \wedge dq$ in canonical Darboux coordinates. We use also the multi-index notation. Let us consider the well-known Moyal $\star$-product of two smooth functions $u,v \in C^\infty(\mathbf R^2)$, defined by
$$u \star v = u.v + \sum_{r \geq 1} \frac{1}{r!}(\frac{1}{2i})^r P^r(u,v),$$ where
$$P^r(u,v) := \Lambda^{i_1j_1}\Lambda^{i_2j_2}\dots \Lambda^{i_rj_r}\partial_{i_1i_2\dots i_r} u \partial_{j_1j_2\dots j_r}v,$$ with $$\partial_{i_1i_2\dots i_r} := \frac{\partial^r}{\partial x^{i_1}\dots \partial x^{i_r}}, x:= (p,q) = (p_1,\dots,p_n,q^1,\dots,q^n)$$ as multi-index notation. It is well-known that this series converges in the Schwartz distribution spaces $\mathcal S (\mathbf R^n)$. We apply this to the special case $n=1$. In our case we have only $x = (x^1,x^2) = (p,q)$.
\begin{proposition}\label{3.1}  \Large\bf\sf
In the above mentioned canonical Darboux coordinates $(p,q)$ on the orbit $\Omega_+$, the Moyal $\star$-product satisfies the relation
$$i\tilde{Z} \star i\tilde{T} - i\tilde{T} \star i\tilde{Z} = i\widetilde{[Z,T]}, \forall Z, T \in \aff(\mathbf R).$$
\end{proposition}

Consequently, to each adapted chart $\psi$ in the sense of \cite{arnalcortet2}, we associate a $G$-covariant $\star$-product.

\begin{proposition}[see \cite{gutt}]     \Large\bf\sf
Let $\star$ be a formal differentiable $\star$-product on $C^\infty(M, \mathbf R)$, which is covariant under $G$. Then there exists a representation $\tau$ of $G$ in $\Aut N[[\nu]]$ such that
$$\tau(g)(u \star v) = \tau(g)u \star \tau(g)v.$$
\end{proposition}

Let us denote by $\mathcal F_pu$ the partial Fourier transform of the function $u$ from the variable $p$ to the variable $x$, i.e.
$$\mathcal F_p(u)(x,q) := \frac{1}{\sqrt{2\pi}}\int_{\mathbf R} e^{-ipx} u(p,q)dp.$$ Let us denote by $ \mathcal F^{-1}_p(u) (x,q)$ the inverse Fourier transform. 
\begin{lemma}\label{lem3.1}   \Large\bf\sf
1. $\partial_p \mathcal F^{-1}_p(p.u) = i \mathcal F^{-1}_p(x.u)  $ ,

2. $ \mathcal F_p(v) = i \partial_x\mathcal F_p(v)  $ ,

3. $P^k(\tilde{Z},\mathcal F^{-1}_p(u)) = (-1)^k \beta e^q \frac{\partial^k\mathcal F^{-1}_p(u)}{\partial^kp}, \mbox{ with } k \geq 2.$
\end{lemma}

For each $Z \in \aff(\mathbf R)$, the corresponding Hamiltonian function is $\tilde{Z} =  \alpha p + \beta e^q $ and we can consider the operator $\ell_Z$ acting on dense subspace $L^2(\mathbf R^2, \frac{dpdq}{2\pi})^\infty$ of smooth functions by left $\star$-multiplication by $i \tilde{Z}$, i.e. $\ell_Z(u) = i\tilde{Z} \star u$. It is then continued to the whole space $L^2(\mathbf R^2, \frac{dpdq}{2\pi})$. It is easy to see that, because of the relation in Proposition (\ref{3.1}), the correspondence $Z \in \aff(\mathbf R) \mapsto \ell_Z = i\tilde{Z} \star .$ is a representation of the Lie algebra $\aff(\mathbf R)$ on the space $N[[\frac{i}{2}]]$ of formal power series in the parameter $\nu = \frac{i}{2}$ with coefficients in $N = C^\infty(M,\mathbf R)$, see e.g. \cite{gutt} for more detail.

We study now the convergence of the formal power series. In order to do this, we look at the $\star$-product of $i\tilde{Z}$ as the $\star$-product of symbols and define the differential operators corresponding to $i\tilde{Z}$. It is easy to see that the resulting correspondence is a representation of $\mathfrak g $ by pseudo-differential operators. 

\begin{proposition}  \Large\bf\sf
For each $Z \in \aff(\mathbf R)$ and for each compactly supported $C^\infty$ function $u \in C^\infty_0(\mathbf R^2)$, we have
$$\hat{\ell}_Z(u) := \mathcal F_p \circ \ell_Z \circ \mathcal F^{-1}_p(u) = \alpha (\frac{1}{2}\partial_q - \partial_x)u + i\beta e^{q -\frac{x}{2}}u.$$
\end{proposition}

{\it The associate irreducible unitary representations}

Our aim in this section is to exponentiate the obtained representation $\hat{\ell}_Z$ of the Lie algebra $\aff(\mathbf R)$ to the corresponding representation of the Lie group $\Aff_0(\mathbf R)$. We shall prove that the result is exactly the irreducible unitary representation $T_{\Omega_+}$ obtained from the orbit method or Mackey small subgroup method applied to this group $\Aff(\mathbf R)$.
Let us recall first the well-known list of all the irreducible unitary representations of the group of affine transformation of the real straight line.
\begin{theorem} [\cite{gelfandnaimark}]\label{4.1}   \Large\bf\sf
Every irreducible unitary representation of the group $\Aff(\mathbf R)$ of all the affine transformations of the real straight line, up to unitary equivalence, is equivalent to one of the pairwise non-equivalent representations:
\begin{itemize}
\item the infinite dimensional representation $S$, realized in the space $L^2(\mathbf R^*, \frac{dy}{\vert y\vert})$, where $\mathbf R^* = \mathbf R \setminus \{0\}$ and is defined by the formula
$$(S(g)f)(y) := e^{iby}f(ay), \mbox{ where } g = \left(\begin{array}{cc} a & b\\ 0 & 1 \end{array}\right),$$
\item the representation $U^\varepsilon_\lambda$, where $\varepsilon = 0,1$, $\lambda \in \mathbf R$, realized in the 1-dimensional Hilbert space $\mathbf C^1$ and is given by the formula
$$U^\varepsilon_\lambda(g) = \vert a \vert^{i\lambda}(\sgn a)^\varepsilon .$$
\end{itemize}
\end{theorem}
Let us consider now the connected component $G= \Aff_0(\mathbf R)$. The irreducible unitary representations can be obtained easily from the orbit method machinery.
\begin{theorem}   \Large\bf\sf
The representation $\exp(\hat{\ell}_Z)$ of the group $G=\Aff_0(\mathbf R)$ is exactly the irreducible unitary representation $T_{\Omega_+}$ of $G=\Aff_0(\mathbf R)$ associated following the orbit method construction, to the orbit $\Omega_+$, which is the upper half-plane $\mathbf H \cong \mathbf R \rtimes \mathbf R^*$, i. e.
+$$(\exp(\hat{\ell}_Z)f)(y) = (T_{\Omega_+}(g)f)(y) = e^{iby}f(ay),\forall f\in L^2(\mathbf R^*, \frac{dy}{\vert y\vert}), $$ where  $g = \exp Z = \left(\begin{array}{cc} a & b \\ 0 & 1 \end{array}\right).$
\end{theorem}

By analogy, we have also
\begin{theorem}    \Large\bf\sf
The representation $\exp(\hat{\ell}_Z)$ of the group $G=\Aff_0(\mathbf R)$ is exactly the irreducible unitary representation $T_{\Omega_-}$ of $G=\Aff_0(\mathbf R)$ associated following the orbit method construction, to the orbit $\Omega_-$, which is the lower half-plane $\mathbf H \cong \mathbf R \rtimes \mathbf R^*$, i. e.
$$(\exp(\hat{\ell}_Z)f)(y) = (T_{\Omega_-}(g)f)(y) = e^{iby}f(ay),\forall f\in L^2(\mathbf R^*, \frac{dy}{\vert y\vert}), $$ where  $g = \exp Z = \left(\begin{array}{cc} a & b \\ 0 & 1 \end{array}\right).$
\end{theorem}

\subsection{The group of affine transformations of the complex straight line}
Recall that the Lie algebra ${\mathfrak g} = \aff({\bf C})$ of affine transformations of the complex straight line is described as follows, see [D].

It is well-known that the group $\Aff({\bf C})$ is a four (real) dimensional Lie group which is isomorphism to the group of matrices:
$$\Aff({\bf C}) \cong  \left\{\left(\begin{array}{cc} a & b \\ 0 & 1 \end{array}\right) \vert a,b \in {\bf C}, a \ne 0 \right\}$$

The most easy method is to consider $X$,$Y$ as complex generators,
$X=X_1+iX_2$ and $Y=Y_1+iY_2$. Then from the relation $[X,Y]=Y$, we get$ [X_1,Y_1]-[X_2,Y_2]+i([X_1Y_2]+[X_2,Y_1]) = Y_1+iY_2$.
This mean that the Lie algebra $\aff({\bf C})$ is a real 4-dimensional Lie algebra, having 4 generators with the only nonzero Lie brackets: $[X_1,Y_1] - [X_2,Y_2]=Y_1$; $[X_2,Y_1] + [X_1,Y_2] = Y_2$ and we can choose another basic noted again by the same letters to have more clear Lie brackets of this Lie algebra:
$$[X_1,Y_1] = Y_1; [X_1,Y_2] = Y_2; [X_2,Y_1] = Y_2; [X_2,Y_2] = -Y_1$$

\begin{remark}{\Large\bf\sf
Let us denote:
$$\mathbf H_{k} = \{w=q_{1}+iq_{2} \in {\bf C} \vert -\infty< q_1<+\infty ; 2k\pi < q_2< 2k\pi+2\pi\}; k=0,\pm1,\dots$$
$$L=\{{\rho}e^{i\varphi} \in {\bf C} \vert 0< \rho < +\infty; \varphi = 0\} \mbox{ and } {\bf C} _{k } = {\bf C} \backslash L$$
is a univalent sheet of the Riemann surface of the complex variable multi-valued analytic function $\Ln(w)$, ($k=0,\pm 1,\dots$)
Then there is a natural diffeomorphism $w \in \mathbf H_{k} \longmapsto e^{w} \in {\bf C}_k$ with each $k=0,\pm1,\dots.$ Now consider the map:
$${\bf C} \times {\bf C} \longrightarrow \Omega_F = {\bf C} \times {\bf C}^*$$
$$(z,w) \longmapsto (z,e^w),$$
with a fixed $k \in \mathbf Z$. We have a local diffeomorphism
$$\varphi_k: {\bf C} \times {\bf H}_k \longrightarrow {\bf C} \times {\bf C}_k$$
 $$(z,w) \longmapsto (z,e^w) $$
This diffeomorphism $\varphi_k$ will be needed in the all sequel.
}\end{remark}

On ${\bf C}\times {\bf H}_k$ we have the natural symplectic form
\begin{equation}\omega = \frac{1 }{2}[dz \wedge dw+d\overline {z} \wedge d\overline {w}],\end{equation} induced from $\mathbf C^2$.
Put $z=p_1+ip_2,w=q_1+iq_2$ and $(x^1,x^2,x^3,x^4)=(p_1,q_1,p_2,q_2) \in {\bf R}^4$, then
$$\omega = dp_1 \wedge dq_1-dp_2 \wedge dq_2.$$ The corresponding symplectic matrix of $\omega$ is
$$ \wedge = \left(\begin{array}{cccc} 0 & -1 & 0 & 0 \\
		           1 & 0 & 0 & 0 \\
                                                0 & 0 & 0 & 1 \\
                                                0 & 0 & -1 & 0 \end{array}\right)
\mbox{   and   }
 \wedge^{-1} = \left(\begin{array}{cccc} 0 & 1 & 0& 0 \\
		           -1 & 0 & 0 & 0 \\
                                                0 & 0 & 0 & -1 \\
                                                0 & 0 & 1 & 0 \end{array}\right)$$

We have therefore the Poisson brackets of functions as follows. With each $f,g \in {\bf C}^{\infty}(\Omega)$
$$\{f,g\} = \wedge^{ij}\frac{\partial f }{\partial x^i}\frac{\partial g}{\partial x^j} =  
\wedge^{12}\frac{\partial f }{\partial p_1}\frac{\partial g}{\partial q_1}+
\wedge^{21}\frac{\partial f }{\partial q_1}\frac{\partial g}{\partial p_1} +
\wedge^{34}\frac{\partial f}{\partial p_2}\frac{\partial g}{\partial q_2} +
\wedge^{43}\frac{\partial f}{\partial q_2}\frac{\partial g}{\partial p_2} = $$
$$\ \ \ \ \ \ \ =\frac{\partial f }{\partial p_1}\frac{\partial g}{\partial q_1} -
\frac{\partial f}{\partial q_1}\frac{\partial g}{\partial p_1} -
\frac{\partial f}{\partial p_2}\frac{\partial g }{\partial q_2} +
\frac{\partial f}{\partial q_2}\frac{\partial g }{\partial p_2} = $$
$$\ \ =2\Bigl[\frac{\partial f}{\partial z}\frac{\partial g}{\partial w} -
\frac{\partial f}{\partial w}\frac{\partial g}{\partial z} +
\frac{\partial f}{\partial \overline z}\frac{\partial g}{\partial \overline{w}} -
\frac{\partial f}{\partial \overline w}\frac{\partial g}{\partial \overline z}\Bigr]$$

\begin{proposition} \Large\bf\sf
Fixing the local  diffeomorphism $\varphi_k (k \in {\bf Z})$, we have: 
\begin{enumerate}
\item For any element $A \in \aff(\mathbf C)$, the corresponding Hamiltonian function $\widetilde{A}$  in local coordinates $(z,w)$ of the orbit $\Omega_F$  is of the form
$$\widetilde A\circ\varphi_k(z,w) = \frac{1}{2} [\alpha z +\beta e^w + \overline{\alpha} \overline{z} + \overline{\beta}e^{\overline {w}}]$$
\item In local coordinates $(z,w)$ of the orbit $\Omega_F$, the symplectic Kirillov form $\omega_F$ is just the standard form (1).
\end{enumerate}
\end{proposition}

{\it Computation of Operators $\hat{\ell}_A^{(k)}$.}

\begin{proposition}\label{Proposition 3.1}  \Large\bf\sf
With $A,B \in \aff({\bf C})$, the Moyal $\star$-product satisfies the relation:
\begin{equation} i \widetilde{A} \star i \widetilde{B} - i \widetilde{B} \star i \widetilde{A} = i[\widetilde{A,B} ]\end{equation}
\end{proposition}

For each $A \in \hbox{aff}({\bf C}$), the corresponding Hamiltonian function is
$$\widetilde{A} = \frac{1}{2} [\alpha  z + \beta e^{w} + \overline{\alpha}\overline z + \overline{\beta} e^ {\overline w}] $$ 
and we can consider the operator  ${\ell}_A^{(k)}$  acting on dense subspace
$L^2({\bf R}^2\times ({\bf R}^2)^*,\frac{dp_1dq_1dp_2dq_2}{(2\pi)^2} )^{\infty}$
 of smooth functions by left $\star$-multiplication by $i \widetilde{A}$, i.e:
${\ell}_A^{(k)} (f) = i \widetilde{A} \star f$. Because of the relation in Proposition 3.1, we have
\begin{corollary}\label{Consequence 3.2} \Large\bf\sf
 \begin{equation}{\ell}_{[A,B]}^{(k)} = {\ell}_A^{(k)} \star {\ell}_B^{(k)} - {\ell}_B^{(k)} \star {\ell}_A^{(k)} := {\Bigl[ {\ell}_A^{(k)},  {\ell}_B^{(k)}\Bigr]}^{\star}\end{equation}
\end{corollary}

From this it is easy to see that, the correspondence $A \in \aff({\bf C}) \longmapsto {\ell}_A^{(k)} = $i$\widetilde {A} \star$. is a representation of the Lie algebra $\aff({\bf C}$) on the space N$\bigl[[\frac{i}{2}]\bigr] $ of formal power series, see [G] for more detail.

\begin{proposition}\label{Proposition 3.4}   \Large\bf\sf
For each $A = \left(\begin{array}{cc}\alpha & \beta \cr 0 & 0 \cr\end{array}\right) \in \aff({\bf C}) $
 and for each compactly supported $C^{\infty}$-function $f \in C_0^{\infty}({\bf C} \times {\bf H}_k)$, we have:
\begin{equation} {\ell}_A^{(k)}{f} := {\mathcal F}_z \circ \ell_A^{(k)} \circ {\mathcal F}_z^{-1}(f) = [\alpha (\frac{1}{2} \partial_w - \partial_{\overline \xi})f + \overline \alpha (\frac{1 }{2}\partial_{\overline w} - \partial_\xi)f + \end{equation}
$$+\frac{i}{2}(\beta e^{w-\frac{1}{2}\overline \xi} + \overline \beta e^{\overline w - \frac{1}{2} \xi})f] $$
\end{proposition}

\begin{remark}\label{Remark 3.5}{\Large\bf\sf  Setting new variables  u = $w - \frac{1}{ 2}\overline{\xi}$;$v = w + \frac{1 }{2}{\overline{\xi}}$ we have
\begin{equation}\hat {\ell}_A^{(k)}(f) = \alpha\frac{ \partial f }{\partial u}+ \overline{\alpha}\frac{\partial f }{\partial{\overline{u}}}+ \frac{i }{2}(\beta e^{u}+\overline{\beta}e^{\overline{u}})f \vert_{(u,v)}\end{equation}
i.e $\hat {\ell}_A^{(k)} = \alpha\frac{ \partial }{\partial u}+ \overline{\alpha}\frac{ \partial  }{\partial{\overline{u}}}+ \frac{i }{2}(\beta e^{u}+\overline{\beta}e^{\overline{u}})$,which provides a ( local) representation of the Lie algebra  aff({\bf C}).
}\end{remark}

{\it The Irreducible Representations of $\widetilde{\Aff}({\bf C})$. }
Since $\hat {\ell}_A^{(k)}$ is a representation of the Lie algebra  $\widetilde{\hbox {Aff}} ({\bf C})$, we have:
$$\exp(\hat {\ell}_A^{(k)}) = \exp\bigl(\alpha\frac{ \partial }{\partial {u}}+ \overline{\alpha}\frac{ \partial  }{\partial{\overline{u}}}+ \frac{i }{2}(\beta e^{u}+\overline{\beta}e^{\overline{u}})\bigr)$$ is just the corresponding representation of the corresponding connected and simply connected Lie group $\widetilde\Aff ({\bf C})$.

Let us first recall the well-known list of all the irreducible unitary representations of the group of affine transformation of the complex straight line, see [D] for more details.

\begin{theorem}\label{Theorem 4.1} \Large\bf\sf
Up to unitary equivalence, every irreducible unitary representation of $\widetilde{\hbox {Aff}} ({\bf C})$ is unitarily equivalent to one the following one-to-another non-equivalent irreducible unitary representations:
\begin{enumerate}
\item The unitary characters of the group, i.e the one dimensional unitary representation $U_{\lambda},\lambda \in {\bf C}$, acting in ${\bf C}$ following the formula
$U_{\lambda}(z,w) = e^{{i\Re(z\overline{\lambda})}}, \forall (z,w) \in \widetilde{\Aff} ({\bf C}), \lambda \in {\bf C}.$
\item The infinite dimensional irreducible representations $T_{\theta},\theta \in {\mathbf S}^1$, acting on the Hilbert space $L^{2}(\mathbf R\times \mathbf S ^1)$ following the formula:
\begin{equation}\Bigr[T_{\theta}(z,w)f\Bigl](x) = \exp \Bigr(i(\Re(wx)+2\pi\theta[\frac{\Im(x+z) }{2\pi}])\Bigl)f(x\oplus z),\end{equation}
Where \ $(z,w) \in\widetilde{\Aff}({\bf C})$  ;  $x \in {\bf R}\times {\mathbf S} ^1= {\bf C} \backslash \{0\}; f \in L^{2}({\bf R}\times {\mathbf S} ^1);$
$$ x\oplus z = Re(x+z) +2 \pi i \{\frac{\Im(x+z) }{2\pi}\}$$
\end{enumerate}
\end{theorem}
In this section we will prove the following important Theorem which
is very interesting for us both in theory and practice.
\par
\begin{theorem}\label{Theorem 4.2} \Large\bf\sf
The representation $\exp(\hat {\ell}_A^{(k)})$ of the group $\widetilde{\Aff}({\bf C})$ is the irreducible unitary representation
$T_\theta$ of $\widetilde{\Aff}({\bf C})$ associated, following the orbit method construction, to the orbit $\Omega$, i.e:
$$\exp(\hat {\ell}_A^{(k)})f(x) = [T_\theta (\exp A)f](x),$$
where $f \in L^{2}({\bf R}\times {\mathbf S} ^1) ; A = \left(\begin{array}{cc}\alpha & \beta \cr 0 & 0 \cr\end{array}\right) \in \aff({\bf C}) ; \theta \in {\mathbf S}^1 ; k = 0, \pm1,\dots$
\end{theorem}

\begin{remark}\label{Remark 4.3} {\Large\bf\sf
We say that a real Lie algebra ${\mathfrak g}$ is in the class $\overline{MD}$ if every K-orbit is of dimension, equal 0 or dim ${\mathfrak g}$. Further more, one proved that
([D, Theorem 4.4]) 
Up to isomorphism, every Lie algebra of class $\overline {MD}$ is one of the following:
\begin{enumerate}
\item Commutative Lie algebras.
\item Lie algebra $\aff({\bf R})$ of affine transformations of the real straight line
\item Lie algebra $\aff({\bf C})$ of affine transformations of the complex straight line.
\end{enumerate}
Thus, by calculation for the group of affine transformations of the real straight line $\Aff({\bf R})$ in [DH] and here for the group affine transformations of the complex straight line $\Aff({\bf C})$ we obtained  a description of the quantum $\overline {MD}$ co-adjoint orbits.
}\end{remark}

\subsection{$MD_4$-groups}
We refer the reader to the results of Nguyen Viet Hai \cite{hai3}-\cite{hai4} for the class of $MD_4$-groups (i.e. 4-dimensional solvable Lie groups, all the co-adjoint of which are of dimension 0 or maximal). It is interesting that here he obtained the same exact computation for $\star$-products and all representations.

\subsection{$SO(3)$} As an typical example of compact Lie group, the author proposed Job A. Nable to consider the case of $SO(3)$. We refer the reader to the  results of Job Nable \cite{nable1}-\cite{nable3}. In these examples, it is interesting that the $\star$-products, in some how as explained in these papers, involved the Maslov indices and Monodromy Theorem.

	\subsection{Exponential groups}
Arnal-Cortet constructed star-products for this case \cite{arnalcortet1}-\cite{arnalcortet2}.
	\subsection{Compact groups}
We refer readers to the works of C. Moreno \cite{moreno}.

\section{Algebraic Noncommutative Chern Characters}

Let $G$ be a compact group, $\HP_*(C^*(G))$ the periodic cyclic homology
introduced
in \S2. Since $C^*(G) = \lim_{\rightarrow\atop N} \prod_{i=1}^N
\Mat_{n_i}({\mathbf C})$, $\HP_*(C^*(G))$ coincides with the $\HP_*(C^*(G))$
defined by J. Cuntz-D. Quillen
\cite{CQ}.

\begin{lem}\label{31}    \Large\bf\sf
Let $\{I_N\}_{N\in {\mathbf N}}$ be the above defined collection of ideals
in $C^*(G)$. Then $$\K_*(C^*(G)) = \lim_{\rightarrow \atop N\in {\mathbf N}}
\K_*(I_N) = \K_*({\mathbf C}(\mathbf T)),$$ where $\mathbf T$ is the fixed
maximal torus in $G$. \end{lem}  First note that the algebraic K-theory of
C*-algebras has the stability property $$\K_*(A \otimes M_n(\mathbf C)) \cong
\K_*(\mathbf C(\mathbf T)).$$ Hence, $$\lim_{\rightarrow} \K_*(I_{n_i}) \cong
\K_*(\prod_{w=\mbox{ highest weight}} {\mathbf C}_w)\cong \K_*(\mathbf C(\mathbf
T)),$$ by Pontryagin duality.

J. Cuntz and D. Quillen \cite{CQ} defined the so called $X$-complexes of
${\mathbf C}$-algebras and then used some ideas of Fedosov product to
define algebraic Chern
characters. We now briefly recall the their definitions. For a
(non-commutative) associate ${\mathbf C}$-algebra $A$, consider the space
of even
non-commutative differential forms $\Omega^+(A) \cong RA$, equipped with
the Fedosov
product
$$\omega_1 \circ \omega_2 := \omega_1\omega_2 - (-1)^{\vert \omega_1\vert}
d\omega_1
d\omega_2,$$ see \cite{CQ}. Consider also the ideal $IA := \oplus_{k\geq 1}
\Omega^{2k}(A)$. It is easy to see that $RA/IA \cong A$ and that $RA$
admits the
universal property that any based linear map $\rho : A \to M$ can be
uniquely extended to a
derivation $D : RA \to M$. The derivations $D : RA \to M$ bijectively
correspond to lifting
homomorphisms from $RA$ to the semi-direct product $RA \oplus M$, which also
bijectively correspond to linear map $\bar\rho : \bar{A}= A /{\mathbf C}
\to M$ given by $$
a\in \bar{A} \mapsto D(\rho a).$$ From the universal property of
$\Omega^1(RA)$, we
obtain a bimodule isomorphism $$RA \otimes \bar{A} \otimes RA \cong
\Omega^1(RA).$$ As in \cite{CQ}, let $\Omega^-A = \oplus_{k\geq 0}
\Omega^{2k+1}A$. Then we have $$\Omega^{-}A \cong RA \otimes \bar{A} \cong
\Omega^1(RA)_\# := \Omega^1(RA)/[(\Omega^1(RA),RA)].$$
\par
J. Cuntz and D. Quillen proved
\begin{thm}(\cite{CQ}, Theorem1): \Large\bf\sf
There exists an isomorphism of
${\mathbf Z}/(2)$-graded complexes
$$\Phi : \Omega A = \Omega^+A \oplus \Omega^{-}A \cong RA \oplus
\Omega^1(RA)_\#,$$ such that
$$\Phi : \Omega^+A \cong RA,$$ is defined by $$\Phi(a_0da_1\dots da_{2n} =
\rho(a_1)\omega(a_1,a_2) \dots \omega(a_{2n-1},a_{2n}),$$
and $$ \Phi : \Omega^{-}A \cong \Omega^1(RA)_\#,$$ $$\Phi(a_0da_1\dots
da_{2n+1})
= \rho(a_1)\omega(a_1,a_2)\dots \omega(a_{2n-1},a_{2n})\delta(a_{2n+1}).$$
With
respect to this identification, the product in $RA$ is just the Fedosov
product on even
differential forms and the differentials on the $X$-complex
$$X(RA) : \qquad RA\cong \Omega^+A \to \Omega^1(RA)_\# \cong \Omega^{-}A
\to RA
$$ become the operators
$$\beta = b - (1+\kappa)d : \Omega^{-}A \to \Omega^+A,$$ $$\delta =
-N_{\kappa^2} b
+ B : \Omega^+A \to \Omega^{-}A,$$ where $N_{\kappa^2} =
\sum_{j=0}^{n-1} \kappa^{2j}$, $\kappa(da_1\dots da_n) := da_n\dots da_1$.
\end{thm}
\par
Let us denote by $IA \triangleleft RA$ the ideal of even non-commutative
differential forms
of
order $\geq 2$. By the universal property of $\Omega^1$ $$\Omega^1(RA/IA) =
\Omega^1RA/((IA)\Omega^1RA + \Omega^1RA.(IA) + dIA).$$ Since $\Omega^1RA =
(RA)dRA = dRA.(RA)$, then $$\Omega^1RA(IA) \cong IA\Omega^1RA \;mod\;
[RA,\Omega^1R].$$
$$\Omega^1(RA/IA)_\# = \Omega^1RA /([RA,\Omega^1RA]+IA.dRA + dIA).$$ For
$IA$-adic tower $RA/(IA)^{n+1}$, we have the complex ${\mathcal
X}(RA/(IA)^{n+1}) :$ $$  RA/IA^{n+1} \leftarrow
\Omega^1RA/([RA,\Omega^1RA]+(IA)^{n+1}dRA + d(IA)^{n+1}).$$
Define
${\mathcal X}^{2n+1}(RA,IA) :$ $$ RA/(IA)^{n+1} \to
\Omega^1RA/([RA,\Omega^1RA]+(IA)^{n+1}dRA + d(IA)^{n+1}) $$ $$\to
RA/(IA)^{n+1},$$
and ${\mathcal X}^{2n}(RA,IA):$ $$ RA/((IA)^{n+1} +[RA,IA^n]) \to
\Omega^1RA/([RA,\Omega^1RA]+d(IA)^ndRA)$$ $$\to RA/((IA)^{n+1}
+[RA,IA^n]).$$
One has
$$b((IA)^ndIA) = [(IA)^n,IA] \subset (IA)^{n+1},$$ $$d(IA)^{n+1} \subset
\sum_{j=0}^n (IA)^jd(IA)(IA)^{n-j} \subset (IA)^n dIA
+ [RA,\Omega^1RA].$$
and hence
$${\mathcal X}^1(RA,IA = X(RA,IA),$$
$${\mathcal X}^0(RA,IA) = (RA/IA)_\#.$$
There is an sequence of maps between complexes $$\dots \to X(RA/IA) \to
{\mathcal
X}^{2n+1}(RA,IA)\to {\mathcal X}^{2n}(RA,IA) \to X(RA/IA) \to \dots $$ We
have the
inverse limits
$$\hat{X}(RA,IA) := \lim_\leftarrow X(RA/(IA)^{n+1}) = \lim_\leftarrow {\mathcal
X}^n(RA,IA).$$
Remark that
$${\mathcal X}^q = \Omega A/F^q\Omega A,$$ $$\hat{X}(RA/IA)= \hat{\Omega}A.$$

We quote the second main result of J. Cuntz and D. Quillen (\cite{CQ},
Thm2), namely:

$$H_i\hat{\mathcal X}(RA,IA) = \HP_i(A).$$

We now apply  this machinery to our case. First we have the following.
\begin{lem}      \Large\bf\sf
$$\lim_{\rightarrow\atop N} \HP^*(I_N) \cong \HP^*({\mathbf C}({\mathbf T})).$$ 
\end{lem}  By similar arguments as in the previous lemma \ref{31}. More
precisely, we
have $$\HP(I_{n_i}) = \HP(\prod_{w=\mbox{highest weight}} {\mathbf C}_w)
\cong
\HP({\mathbf C}({\mathbf T}))$$ by Pontryagin duality.

Now, for each idempotent $e\in M_n(A)$ there is an unique element $x\in
M_n(\widehat{RA})$.
Then the element $$\tilde{e} := x + (x-\frac{1}{2})\sum_{n\geq 1} \frac{2^n(2n-
1)!!}{n!}(x-x^{2n})^{2n}\in M_n(\widehat{RA})$$ is a lifting of $e$ to an
idempotent
matrix in $M_n(\widehat{RA})$. Then the map $[e] \mapsto tr(\tilde{e})$
defines the map
$\K_0 \to H_0(X(\widehat{RA})) = \HP_0(A)$. To an element $g\in\GL_n(A)$ one
associates an element $p\in \GL(\widehat{RA})$ and to the element $g^{-1}$
an element
$q\in \GL_n(\widehat{RA})$ then put
$$x = 1- qp, \mbox{ and } y = 1-pq.$$
And finally, to each class
$[g]\in \GL_n(A)$ one associates $$tr(g^{-1}dg) = tr(1-x)^{-1}d(1-x) =
d(tr(log(1-x))) =  $$ $$ =
-tr\sum_{n=0}^\infty x^ndx\in \Omega^1(A)_\#.$$ Then $[g] \to tr(g^{-1}dg)$
defines
the map $\K_1(A) \to \HH_1(A) = H_1(X(\widehat{RA})) = \HP_1(A)$.

Let $\HP(I_{n_i})$ be the periodic cyclic cohomology defined by
Cuntz-Quillen. Then the pairing
$$\K_*^{alg}(C^*(G)) \times \bigcup_N \HP^*(I_N) \to {\mathbf C}$$ defines an
algebraic non-commutative Chern character $$ch_{alg} : \K_*^{alg}(C^*(G)) \to
\HP_*(C^*(G)),$$ which gives us a variant of non-commutative Chern
characters with
values in $\HP$-groups.

\begin{thm}   \Large\bf\sf
Let $G$ be a compact group and ${\mathbf T}$ a fixed maximal compact torus
of $G$.
Then, the Chern character $$ch_{alg} : \K_*(C^*(G)) \to
\HP_*(C^*(G))$$ is an isomorphism, which can be identified with the
classical Chern
character $$ch: \K_*({\mathbf C}({\mathbf T})) \to \HP_*({\mathbf C}({\mathbf
T}))$$
which is also an isomorphism.
\end{thm}

{

\end{document}